\documentclass[11pt]{amsart}
\usepackage{amssymb}
\usepackage{epsfig}
\usepackage{graphicx}
\newtheorem{theorem}{Theorem}

\newtheorem{corollary}[theorem]{Corollary}

\begin{document}

\title{On colored quandle longitudes and its applications to tangle embeddings and virtual knots}

\author{Maciej Niebrzydowski}
\address{The George Washington University,
Washington, DC 20052}
\email{mniebrz@gmail.com}

\date{May 16th, 2006}

\dedicatory{Dedicated to Professor Louis Kauffman on the occasion
     of his $60$th birthday}

\thispagestyle{empty}

\begin{abstract}
Given a long knot diagram $D$ and a finite quandle $Q$, we consider the set of all quandle colorings of $D$ with a fixed color $q$ of its initial arc. Using this set we define the family 
$\Phi_{Q} ^{q}(K)$ of quandle automorphisms which is a knot invariant. 
For every element $x\in Q$ one can consider the formal sum 
$S^{x}_{\Phi}(K)=\sum_{\phi} \phi(x)$, taken over all $\phi \in \Phi_{Q}^{q}$. 
Such formal sums can be applied to a tangle embedding problem and recognizing non-classical virtual knots.
\end{abstract}

\maketitle

\section{Introduction}
Let $\pi_{K}$ denote the fundamental group of the complement of a knot $K$ in $\mathbb{S}^3$, $m_{K}$ be the meridian of $K$ and $l_{K}$ its longitude (see \cite{B-Z} for definitions).
It is well known that, up to isomorphism, the triple $(\pi_{K},\,m_{K},\,l_{K})$ is a classifying invariant of the knot:
\begin{theorem}[\cite{W}]
Two oriented knots $K$ and $K'$ are isotopic if and only if there is an isomorphism $\phi\colon \pi_K\to
\pi_{K'}$ such that $\phi(m_K)=m_{K'}$ and $\phi(l_K)=l_{K'}$.
\end{theorem}
Given the power of the group system $(\pi_{K},\,m_{K},\,l_{K})$, it is natural to try to use the pair $(m_{K},\,l_{K})$ when designing knot invariants.
In \cite{Eis}, M. Eisermann defined the following invariants. For a finite group $G$ and an   element $x\in G$, consider the set of all representations $\rho\colon \pi_K\to G$, with 
$\rho(m_K)=x$. The sum $\sum_{\rho}\rho(l_K)$, which is an element of the group ring $\mathbb{Z} G$, is a knot invariant. These invariants were successfully applied to detecting chirality and non-inversibility of some knots.
 
In this paper, we change this idea so that it can incorporate quandles, and then use it to find obstructions to tangle embeddings. Our invariants can be easily generalized to the virtual category, and give some information about non-classicality of virtual knots.

Now let us recall the definition of a quandle as given in \cite{Joy}.\\
\textbf{Definition.} A \textit{quandle} is a set $Q$ equipped with two binary operations 
$\ast,\,\bar{\ast}$, satisfying the following identities:
\begin{enumerate}
\item [Q1.] $x\ast x=x$, for all $x\in Q$;
\item [Q2.] $(x\ast y)\,\bar{\ast}\,y=x=(x\,\bar{\ast}\,y)\ast y$, for all $x,\,y\in Q$;
\item [Q3.] $(x\ast y)\ast z=(x\ast z)\ast (y\ast z)$, for all $x,\,y,\,z\in Q$.
\end{enumerate}
Let the symbol $\ast ^\epsilon \in \{\ast,\,\bar{\ast}\}$ denote a generic quandle operation, so that the meaning of $x\ast ^\epsilon y$ is either $x\ast y$ or $x\bar{\ast}\,y$. We will use subscripts to specify several possibly different operations. If the symbol $\ast ^\epsilon$ is 
already used, then the symbol $\ast ^{-\epsilon}$ denotes the opposite operation. For example, if
$\ast ^\epsilon=\ast$, then $\ast ^{-\epsilon}=\bar{\ast}$.  
We use a standard convention for products in non-associative algebras, called the left normed convention i.e., whenever parentheses are omitted in a product of elements $a_1$, $a_2,\ldots,$ $a_n$ of $Q$, then $a_1\ast ^{\epsilon_1} a_2\ast ^{\epsilon_2} \ldots 
\ast ^{\epsilon_{n-1}} a_n=((\ldots ((a_1\ast ^{\epsilon_1} a_2)\ast ^{\epsilon_2}
a_3)\ast ^{\epsilon_3} \ldots)\ast ^{\epsilon_{n-2}} a_{n-1})\ast ^{\epsilon_{n-1}} a_n$ (left association) for example, $a\ast b\,\bar{\ast}\,c\ast d=((a\ast b)\,\bar{\ast}\,c)\ast d$.\\
\textbf{Definition.} A \textit{homomorphism of quandles} $Q$ and $Q'$ is a map $f\colon Q\to Q'$ satisfying the condition $f(x\ast y)=f(x)\ast f(y)$, for all $x,\,y\in Q$.
An automorphism of $Q$ is a bijective homomorphism $f\colon Q\to Q$.

A basic and important example of the automorphism of a quandle $Q$ is a function $f_q$ defined by
$f_q(x)=x\ast q$, where $q,\, x\in Q$, or a function $\overline{f_q}$ given by
$\overline{f_q}(x)=x\,\bar{\ast}\,q$.

Initially we will work with classical long knots.\\
\textbf{Definition.} A \textit{long knot diagram} is a smooth immersion $f\colon \mathbb{R}\to \mathbb{R}^2$ with crossing information at each double point, and such that $f(x)=(x,0)$ for
$|x|$ sufficiently large.\\
Long knot diagrams are assumed to be oriented from the left to the right.

\begin{figure}
\begin{center}
\includegraphics[height=2.5 cm]{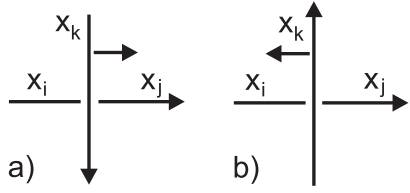}
\caption{Two types of crossings with quandle generators assigned to arcs.\label{cr}}
\end{center}
\end{figure}
Given a diagram $D$ representing the long knot $K$, we can assign to it the \textit{fundamental quandle} $Q(K)$. Its generators correspond to arcs of the diagram, and relations correspond to crossings and are of the form: $x_i\ast x_k=x_j$ if crossing is as in Fig.\ref{cr}(a) or
$x_i\,\bar{\ast}\,x_k=x_j$ in the case of crossing depicted in Fig.\ref{cr}(b).
$Q(K)$ does not change (up to isomorphism) under Reidemeister moves, so it does not depend on the long knot diagram that we choose (see \cite{Joy} for details).

We can obtain a closed knot $\tilde{K}$ from a long knot $K$ by joining the initial and the last arcs of $K$. To define the fundamental quandle of $\tilde{K}$, we only need to add to $Q(K)$ one relation identifying generators corresponding to arcs that are being joined.\\
\begin{figure}[ht]
\begin{center}
\includegraphics[height=2.5 cm]{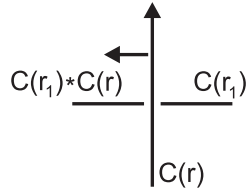}
\caption{Quandle coloring relation.\label{col}}
\end{center}
\end{figure}
The following definition applies to both long and closed knots.\\
\textbf{Definition.} Let $X$ be a fixed quandle and $D$ be a given diagram of an oriented knot $K$. Let $R$ be the set of arcs of this diagram. The normals of arcs are chosen so that the pair (tangent, normal) matches the usual orientation of the plane. A \textit{quandle coloring} $C$ is a map $C\colon R\to X$ such that at every crossing the following relation is satisfied: if $r$ is the over-arc at a crossing, and $r_1$ is the under-arc such that the normal of $r$ points away from it, then the second under-arc should receive color $C(r_1)\ast C(r)$ (see Fig.\ref{col}).
For a survey on quandles and quandle colorings see \cite{C-K-S}.
Alternatively, quandle coloring can be viewed as a homomorphism from the fundamental quandle $Q(K)$ to the quandle $X$.
\section{Colored quandle longitudes}
Let $D$ be a diagram of the long knot $K$ with the initial arc denoted as $r_1$ and $Q$ be a finite quandle. Let us fix an element $q\in Q$. We can consider the set $Col(D,Q,q)$ of all colorings
of $D$ with quandle $Q$, satisfying the condition that the color of $r_1$ is $q$.\\
The cardinality of this set is a knot invariant, but we would like to obtain stronger invariant by means of longitudinal information.

\begin{figure}
\begin{center}
\includegraphics[height=4.5cm]{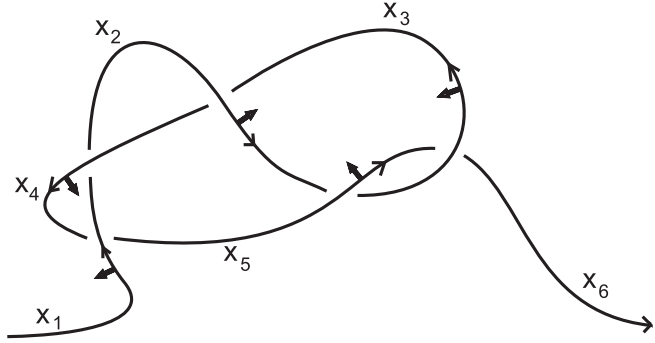}
\caption{A long diagram of the knot $5_2$.\label{k5_2}}
\end{center}
\end{figure}
In the context of the fundamental quandle $Q(K)$, the longitude can be described as follows.\\
\textbf{Definition.}
Travelling along the knot diagram $D$ from the left to the right, we label the arcs of $D$ by generators $x_1,\ldots ,x_{n+1}$ of $Q(K)$, in a consecutive manner (see Fig.\ref{k5_2} for an example). Here, $n$ is the number of crossings of the diagram. Let $x_{o(i)}$ denote the generator assigned to the over-arc that is encountered between arcs labeled by $x_i$ and $x_{i+1}$, and let $\upsilon(i)$ be equal to 1 if the normal to this over-arc points away from arc labeled $x_i$ (Fig.\ref{cr}(a)) 
and -1 otherwise (Fig.\ref{cr}(b)).\\ 
The \textit{quandle longitude} is an automorphism $\eta\colon Q(K)\to Q(K)$ defined by the formula 
$$\eta(x)=x\ast^{-\epsilon_1}x_1\ast^{\epsilon_1}x_{o(1)}\ast^{-\epsilon_2}x_2
\ast^{\epsilon_2}x_{o(2)}\ast^{-\epsilon_3}\ldots\ast^{-\epsilon_n}x_n
\ast^{\epsilon_n}x_{o(n)},$$
where $\ast^{\epsilon_i}=\ast$ if $\upsilon(i)=1$ and $\ast^{\epsilon_i}=\bar{\ast}$ otherwise.

For example, the quandle longitude for the long knot $5_2$ on Fig.\ref{k5_2} can be written as
$$\eta_{5_2}(x)=x\ast x_1\,\bar{\ast}\,x_4\ast x_2\,\bar{\ast}\,x_5\ast x_3\,\bar{\ast}\,x_2
\ast x_4\,\bar{\ast}\,x_1\ast x_5\,\bar{\ast}\,x_3.$$
To simplify the notation, we will describe the quandle longitude as a list of generators, in exactly the same order as in definition of $\eta$, with the following additional information: if generator $x_i$ is to be preceded by the operation $\bar{\ast}$, it will be written as $\overline{x_i}$, otherwise it will be left as $x_i$.\\
Thus, we write the quandle longitude in our example as
$$\{x_1,\,\overline{x_4},\,x_2,\,\overline{x_5},\,x_3,\,\overline{x_2},\,x_4,\,
\overline{x_1},\,x_5,\,\overline{x_3}\}.$$
For a slightly different definition of the quandle longitude, suitable for framed knots, see 
\cite{K1}.

Now we can unify the idea of the quandle longitude with the idea of quandle colorings.\\
\textbf{Definition.} Let $\zeta$ be the quandle coloring of the long knot diagram $D$, using some fixed quandle Q.
Then the \textit{colored quandle longitude} corresponding to $\zeta$ is the automorphism 
$\phi_{\zeta}\colon Q\to Q$, defined by the formula
$$\phi_{\zeta}(x)=x\ast^{-\epsilon_1}\zeta(x_1)\ast^{\epsilon_1}\zeta(x_{o(1)})\ast^{-\epsilon_2}\zeta(x_2)\ast^{\epsilon_2}\zeta(x_{o(2)})\ast^{-\epsilon_3}\ldots\ast^{-\epsilon_n}
\zeta(x_n)\ast^{\epsilon_n}\zeta(x_{o(n)}),$$
where $x\in Q$, and we apply conventions from the previous definition when determining the type of the operation that is used. The colored quandle longitude is not changed by Reidemeister moves.

If $\zeta$ was the coloring of the long knot $5_2$ from the example, then we would write the corresponding colored quandle longitude as
$$\{\zeta(x_1),\,\overline{\zeta(x_4)},\,\zeta(x_2),\,\overline{\zeta(x_5)},\,\zeta(x_3),\,
\overline{\zeta(x_2)},\,\zeta(x_4),\,\overline{\zeta(x_1)},\,\zeta(x_5),\,
\overline{\zeta(x_3)}\}.$$

Now we are ready to define our invariants.\\
\textbf{Definition.} Let $D$, $K$, $Q$, $Col(D,Q,q)$ be as in the beginning of this section.
Then we can define the \textit{colored quandle longitude invariant} $\Phi_{Q}^{q}(K)$ as the family of all colored quandle longitudes corresponding to the set of colorings
$Col(D,Q,q)$ i.e., $\Phi_{Q}^{q}(K)=\{\phi_{\zeta}\,|\,\zeta\in Col(D,Q,q)\}$.

The next natural step is to choose an element $x\in Q$, and consider the formal sum $$S^{x}_{\Phi}(K)=\sum_{\phi} \phi(x),$$ taken over all $\phi\in \Phi_{Q}^q$, which is an element of the free $\mathbb{Z}$-module generated by elements of $Q$.

The above invariants were defined for long knots, but it is well known that in the classical case the long knot theory coincides with the theory of closed knots. Any two long knots obtained from the closed knot by breaking it at two different points are the same, one can be obtained from the other by a sequence of Reidemeister moves. Thus, our invariants are well defined for classical closed knots, for their value is the same for all long knots obtained from the closed knot by breaking it at some point.

Let us apply these invariants to show that the knot $5_2$ is not equivalent to its mirror image.\\
\textbf{Example.} We choose the coloring quandle $Q$ to be the conjugacy class of the element
$q=(1,2)(3,4,5)$ in the symmetric group $S_5$, with quandle operations defined as 
$a\ast b=b^{-1}ab$ and $a\,\bar{\ast}\,b=bab^{-1}$, where on the right hand side we have a group multiplication in $S_5$. $Q$ has 20 elements. As an element $x$ that will be acted upon by colored quandle longitudes, we choose permutation $(1,2,3)(4,5)$. There are 7 colorings of the long knot $5_2$ (Fig.\ref{k5_2}) with $q$ being the color of the initial arc, therefore we have 7 quandle automorphisms acting on $x$.
The value of the invariant is
$$S^{x}_{\Phi}(5_2)=6\cdot (1,2,4)(3,5)+(1,2,3)(4,5).$$
The element $(1,2,3)(4,5)$ in the above sum is a result of acting on $x$ by the colored quandle longitude corresponding to monochromatic coloring (i.e., when all arcs receive color $q$). It is an identity automorphism.\\
Let $\overline{5_2}$ denote the mirror image of the knot $5_2$.
The value of the sum for this knot is
$$S^{x}_{\Phi}(\overline{5_2})=6\cdot (1,2,5)(3,4)+(1,2,3)(4,5).$$
Since these values are different, the knot $5_2$ is chiral.

We can apply similar technique to show chirality of the knot $9_{42}$ that is known to be the smallest chiral knot not distinguished from its mirror image by the Homflypt and Kauffman polynomials.\\
In this case we take $Q$ to be the alternating group $A_5$ with conjugation as the quandle operation. Let $q=(1,2,3)$ and $x=(2,3,4)$. The values of the invariant are as follows:
$$S^{x}_{\Phi}(9_{42})=7\cdot (2,3,4)+6\cdot (1,4,3),$$
$$S^{x}_{\Phi}(\overline{9_{42}})=7\cdot (2,3,4)+6\cdot (1,2,4).$$
These and other computations included in this paper were obtained using the computer algebra system GAP \cite{GAP}.
\section{Applications to tangle embeddings}
A $2n$-tangle consists of $n$ disjoint arcs in the 3-ball. We ask the following question, 
that was first considered by D. Krebes \cite{Kre}. For a given knot $K$ and a tangle 
$T$, can we embed $T$ into $K$ i.e., is there a diagram of $T$ that extends to 
a diagram of $K$?
This problem is important due to its applications in the study of DNA. 
A number of knot invariants have been 
used to find criteria for tangle embeddings (see for example \cite{P-S-W}, \cite{Rub}).

\begin{figure}
\begin{center}
\includegraphics[height=5cm]{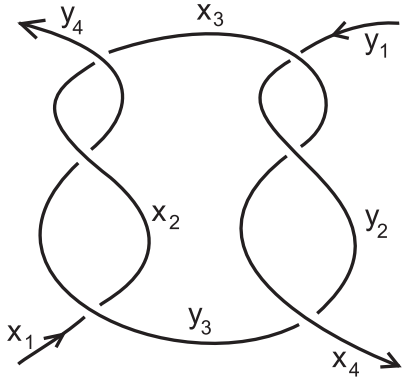}
\caption{Tangle $T_{6_2}$.\label{t62}}
\end{center}
\end{figure}

Here, we use invariants from the previous section to define obstructions to tangle embeddings.
We will work with oriented 4-tangles i.e., tangles that have two inputs and two outputs 
(see Fig.\ref{t62} for an example), but our method is applicable for general $2n$-tangles \footnote{Many authors use the name $n$-tangles for such objects having $n$ inputs and $n$ outputs.}.

First, we need a definition of particular kind of tangle coloring, that is necessary ingredient in what follows.\\
\textbf{Definition.} Let $D_T$ be a tangle diagram, and $Q$ be a quandle. A 
\textit{boundary-monochromatic coloring} of $D_T$ is a map from the set of arcs of $D_T$ to quandle $Q$ satisfying the usual conditions for quandle colorings of knot diagrams, and an additional requirement that all boundary points receive the same color. 

If a tangle $T$ embeds into a knot $K$, then each boundary-monochromatic coloring of $D_T$ can be extended trivially to the whole diagram of $K$ i.e., all arcs outside $D_T$ receive the color of the boundary points of $D_T$. Thus, the existence of nontrivial colorings of $D_T$ gives the first basic obstruction to tangle embeddings, for $T$ can possibly embed only into knots admitting at least the same number of nontrivial colorings (see also \cite{Kre}).

Let us proceed to define our obstructions.\\
Suppose that a tangle $T$ embeds into a long knot $K$ i.e., there exists a 
tangle diagram $D_T$, and a diagram $D$ of $K$ such that $D_T$ is a part of $D$.
Let $Q$ be a finite quandle, $q\in Q$, and let $Col^q_Q(D_T)$ denote the set of all boundary-monochromatic colorings of $D_T$ with colors from $Q$, such that the color of the boundary points of $D_T$ is equal to $q$.\\
Each coloring $\zeta_{D_T}\in Col^q_Q(D_T)$ determines a coloring $\zeta\in Col(D,Q,q)$ (by extending $\zeta_{D_T}$ trivially to $D$).\\
For every such coloring $\zeta$, diagram $D_T$ yields two parts of the colored quandle longitude $\phi_{\zeta}$, obtained by traveling along two components of tangle $T$ according to the orientation of $T$. Let us denote these parts as $\phi_{\zeta}^1$ and $\phi_{\zeta}^2$.

In the case of tangle $T_{6_2}$ from the Fig.\ref{t62} longitude parts would be as follows:
$$\phi_{\zeta}^1=\{\zeta(x_1),\,\overline{\zeta(y_3)},\,\zeta(x_2),\,\overline{\zeta(y_4)},\,
\zeta(x_3),\,\overline{\zeta(y_2)}\},$$ starting from the lower left corner, and
$$\phi_{\zeta}^2=\{\zeta(y_1),\,\overline{\zeta(x_3)},\,\zeta(y_2),\,\overline{\zeta(x_4)},\,
\zeta(y_3),\,\overline{\zeta(x_2)}\},$$ starting from the upper right corner of $D_{T_{6_2}}$.

Notice that we do not know a priori which longitude part appears first in $\phi_{\zeta}$ therefore, we must take into account both possibilities.

Each crossing of $D$ that is outside of $D_T$ contributes $\{q,\,\overline{q}\}$ or 
$\{\overline{q},\, q\}$ to $\phi_{\zeta}$ (all arcs outside $D_T$ have color $q$).
This contribution is trivial because of the second axiom from the definition of quandle.\\
Let $\centerdot$ denote the concatenation of lists. We have the following result.
\begin{theorem}
Let $Q$, $q$, $T$ and $K$ be as above. 
If $T$ embeds into $K$ then either 
$\phi_{\zeta}^1\centerdot \phi_{\zeta}^2$ or $\phi_{\zeta}^2\centerdot \phi_{\zeta}^1$ belongs to $\Phi_{Q}^q$.
\end{theorem}
This observation can be used when utilizing the formal sums $S^{x}_{\Phi}$.
\begin{corollary}
If $T$ embeds into $K$, then either the sum 
$$S^{x}_1(T)=\sum(\phi_{\zeta}^1\centerdot \phi_{\zeta}^2)(x),$$ or the sum
$$S^{x}_2(T)=\sum(\phi_{\zeta}^2\centerdot \phi_{\zeta}^1)(x),$$ taken over all colorings $\zeta$ obtained from the boundary monochromatic colorings of $D_T$, must be included in the formal sum $$S^{x}_{\Phi}(K)=\sum_{\phi} \phi_{\zeta}(x),$$ that is taken over all colorings 
$\zeta\in Col(D,Q,q)$.\\
If this condition is not satisfied for some choice of $Q$, $q$ and $x\in Q$, then tangle $T$ does not embed into the long knot $K$.
\end{corollary}
If the above method gives obstructions to embedding of $T$ into a long knot $K$, then $T$ will also not embed into a closed knot $\tilde{K}$ obtained from $K$ by joining its ends.
Indeed, for any long knot obtained from $\tilde{K}$ by breaking it at some point, the values of the invariants $\Phi_{Q} ^{q}$, and therefore also obstructions to embeddings, are the same.
Finally, if $T$ embeds into $\tilde{K}$, then we can always choose a breaking point that is outside $T$. 

Let us consider an example of this application.\\
\textbf{Example.} As quandle $Q$, we take the alternating group $A_6$ with conjugation as a quandle operation. Let $q=(1,2,3,4)(5,6)$, and $x=(1,2,3,4,5)$.\\
We consider the problem of embedding of tangle $T_{6_2}$ (Fig.\ref{t62}) into a knot $6_3$.
Tangle $T_{6_2}$ has 9 colorings with the color of the boundary points equal to $q$.
The values of the sums $S^{x}_1(T_{6_2})$ and $S^{x}_2(T_{6_2})$ are equal:
$$S^{x}_1(T_{6_2})=S^{x}_2(T_{6_2})=8\cdot (1,2,5,3,4)+(1,2,3,4,5).$$
Any long knot obtained from the closed knot $6_3$ has 33 colorings in $Col(D,Q,q)$. 
All colored quandle longitudes corresponding to these colorings 
act trivially on the element $x$:
$$S^{x}_{\Phi}(6_2)=33\cdot (1,2,3,4,5).$$
Therefore, we conclude that tangle $T_{6_2}$ does not embed into knot $6_3$.
\section{Applications to virtual knots}
\begin{figure}
\begin{center}
\includegraphics[height=4.5cm]{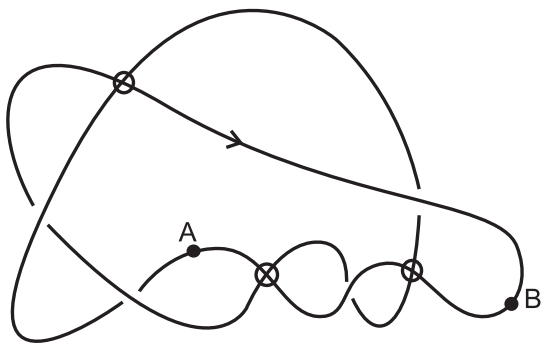}
\caption{A non-classical virtual knot.\label{virt}}
\end{center}
\end{figure}
Virtual knot theory was introduced by Louis Kauffman in \cite{K1}.\\
\textbf{Definition.} A \textit{virtual knot} is defined as an equivalence class of 
4-valent plane diagrams, with an extra crossing information, where a new type of crossing (called virtual crossing and denoted by a small circle around a double point) is allowed.

Virtual knot theory is a generalization of classical knot theory; if two classical knots are equivalent under generalized Reidemeister moves, then they are equivalent under classical ones
(see \cite{G-P-V} for the proof).\\
\textbf{Definition.} By a \textit{long virtual knot diagram} we mean a diagram satisfying the conditions required for a diagram of a classical long knot, but we allow virtual crossings to be present.\\
\textbf{Definition.} A \textit{long virtual knot} is an equivalence class of long virtual
knot diagrams modulo generalized Reidemeister moves (see \cite{Man} for details).

Long virtual knots and their invariants first appeared in \cite{G-P-V}, where it was also shown that breaking the same virtual knot diagram at different points can produce different long virtual knots.

The definition of the fundamental knot quandle, quandle colorings, and invariants that we defined in previous sections can be easily generalized to virtual long knots. We simply ignore the virtual crossings and allow arcs to pass through them.\\
\textbf{Definition.} We say that a virtual knot $\tilde{K}$ is \textit{non-classical} if every diagram of $\tilde{K}$ contains at least one virtual crossing.

Our invariants of long knots can be used to detect virtual knots that are non-classical.
Indeed, it is enough to show that the values of the invariant $S^{x}_{\Phi}$ applied to long
knots obtained by breaking the same virtual knot are not equal.

Another possibility of using these invariants follows from the fact that, unlike in classical case, the connected sum of long virtual knots $K_1$, $K_2$ is not commutative \cite{K05}. If we can show that invariants $\Phi_{Q}^{q}$ of connected sums 
$K_1\# K_2$ and $K_2\# K_1$ are not the same, it will follow that
$K_1$ and $K_2$ are different and both non-classical.\\
\textbf{Example.}Consider the virtual knot $\tilde{K}$ depicted in Fig.\ref{virt}. Let $K_1$ and $K_2$ be long knots obtained from $\tilde{K}$ by breaking it at points $A$ and $B$ respectively.\\
Let $Q$ be the conjugation quandle of the permutation group generated by
$(1,9,6,7,5)(2,10,3,8,4)$ and $(1,10,7,8)(2,9,4,6)$.\\
We choose $q=(2,3,4,5)(6,9,8,10)$, and $x=(2,5)(3,10)(4,9)(7,8)$.\\
Both knots have 17 colorings with the color of the initial arc equal to $q$.
The values of the invariant $S^{x}_{\Phi}$ are as follows:
{\footnotesize
$$S^{x}_{\Phi}(K_1)=(2,6)(3,7)(4,10)(8,9)+(2,9)(4,8)(5,7)(6,10)+2\cdot(1,2)(3,5)(4,7)(6,9)$$$$+2
\cdot(1,2)(3,6)(5,9)(8,10)+2\cdot(1,4)(2,6)(5,10)(7,9)+2\cdot(1,4)(2,9)(3,8)(6,7)$$$$ +(1,7)(2,3)(4,5)(9,10)+(1,7)(2,5)(3,4)(6,8)+(1,8)(3,7)(4,6)(5,9)$$$$+(1,8)(2,10)(3,6)(4,7)
+(1,10)(3,6)(4,9)(5,7)+(1,10)(2,8)(4,7)(5,9)$$$$+(2,5)(3,10)(4,9)(7,8);$$
$$S^{x}_{\Phi}(K_2)=(3,9)(4,7)(5,6)(8,10)+(2,3)(4,6)(5,8)(7,10)+(2,4)(3,5)(6,8)(9,10)$$$$
+(2,6)(3,7)(4,10)(8,9)+(2,7)(3,8)(5,10)(6,9)+(2,8)(3,6)(4,5)(7,9)$$$$ 
+(2,9)(4,8)(5,7)(6,10)+2\cdot(2,10)(3,4)(5,9)(6,7)+(1,6)(2,7)(4,9)(5,8)$$$$
+(1,6)(2,8)(3,10)(5,7)+(1,7)(2,3)(4,5)(9,10)+(1,8)(3,7)(4,6)(5,9)$$$$ 
+(1,8)(2,10)(3,6)(4,7)+(1,9)(4,5)(6,8)(7,10)+(1,10)(2,3)(6,8)(7,9)$$$$
+(2,5)(3,10)(4,9)(7,8).$$}
Since they are different, $\tilde{K}$ is non-classical.

Generalizations of the presented invariants to other quandle-like structures, for example to biquandles, will be investigated in the future paper.


\begin{thebibliography}{99}
\bibitem[B-Z]{B-Z} G. Burde, H. Zieschang, Knots, De Gruyter (1985).
\bibitem[C-K-S]{C-K-S} S. Carter, S. Kamada, M. Saito, Surfaces in 4-space, Springer (2004).
\bibitem[Eis]{Eis} M. Eisermann, Knot colouring polynomials, preprint (2005).
\bibitem[GAP4]{GAP} 
The GAP Group, GAP -- Groups, Algorithms, and Programming, Version 4.4; 2005. (http://www.gap-system.org)  
\bibitem[G-P-V]{G-P-V} M. Goussarov, M. Polyak, O. Viro, Finite type invariants of classical and virtual knots, {\it Topology} \textbf{39} (2000), 1045-1068;\\
e-print: http://arxiv.org/abs/math.GT/9810073
\bibitem[Joy]{Joy} D. Joyce, A classifying invariant of knots, the knot quandle,
{\it Journal of Pure and Applied Algebra} \textbf{23} (1982), 37-65. 
\bibitem[K99]{K1} L. H. Kauffman, Virtual knot theory, {\it European J. Comb.} \textbf{20} (1999), 663-690;\\
e-print: http://arxiv.org/abs/math.GT/9811028
\bibitem[K05]{K05} L. H. Kauffman, Knot diagrammatics, Handbook of knot theory, edited by
W. Menasco and M. Thistlethwaite, Elsevier (2005).
\bibitem[Kre]{Kre} D. A. Krebes, An obstruction to embedding 4-tangles in links,
{\it Journal of Knot Theory and Its Ramifications} \textbf{8} (1999), 321-352;\\
e-print: http://arxiv.org/abs/math.GT/9902119 
\bibitem[Man]{Man} V. O. Manturov, Long virtual knots and their invariants, {\it Journal of Knot Theory and Its Ramifications} \textbf{13} (2004), 1029-1039.
\bibitem[P-S-W]{P-S-W} J. H. Przytycki, D. S. Silver, S. G. Williams, 3-manifolds, tangles and persistent invariants, {\it Math. Proc. Cambridge Phil. Soc.} \textbf{139} (2005), 291-306;\\
e-print: http://arxiv.org/abs/math.GT/0405465  
\bibitem[Rub]{Rub} D. Ruberman, Embedding tangles in links,
{\it Journal of Knot Theory and Its Ramifications} \textbf{9} (2000), 523-530;\\
e-print: http://arxiv.org/abs/math.GT/0001141 
\bibitem[Wal]{W} F. Waldhausen, On irreducible 3-manifolds which are sufficiently large, {\it Ann. of Math.} \textbf{87} (1968), 56-88.
\end{thebibliography}
\end{document}